\input amstex
\input amsppt.sty
\input texdraw

\magnification1200
\hsize14cm
\vsize19cm

\TagsOnRight

\def\StemAE{23}
\def\SlatAC{22}
\def\PoorAB{21}
\def\PeWZAA{20}
\def\PeWiAB{19}
\def\LindAA{18}
\def\KupeAD{17}
\def\KrZeAA{16}
\def\KratBN{15}
\def\KratBH{14}
\def\KratBI{13}
\def\KratBD{12}
\def\KratBG{11}
\def\GrKPAA{10}
\def\GospAB{9}
\def\GeViAA{8}
\def\FuKrAD{7}
\def\FuKrAC{6}
\def\EisTAC{5}
\def\CiKrAD{4}
\def\CiKrAA{3}
\def\CiEKAA{2}
\def\AnBuAA{1}


\def\ZAa{1}
\def\ZA{2}
\def\ZB{3}
\def\ZC{4}
\def\ZD{5}
\def\ZE{6}
\def\ZEa{7}
\def\ZF{8}
\def\ZG{9}
\def\ZH{10}
\def\ZI{11}
\def\ZJ{12}
\def\ZK{13}
\def\ZL{14}
\def\ZM{15}
\def\ZN{16}
\def\ZO{17}
\def\ZP{18}

\def\ZR{20}
\def\ZS{21}
\def\ZT{22}
\def\ZU{23}
\def\ZV{24}
\def\ZW{25}
\def\ZX{26}
\def\ZXa{27}
\def\ZY{28}
\def\ZZ{29}
\def\Za{20}
\def\Zb{31}

\def\TA{1}

\def\po#1#2{(#1)_#2}
\def\({\left(}
\def\){\right)}
\def\[{\left[}
\def\]{\right]}
\def\al{\alpha}
\def\be{\beta}
\def\fl#1{\left\lfloor#1\right\rfloor}

\catcode`\@=11
\font\tenln    = line10
\font\tenlnw   = linew10

\newskip\Einheit \Einheit=0.5cm
\newcount\xcoord \newcount\ycoord
\newdimen\xdim \newdimen\ydim \newdimen\PfadD@cke \newdimen\Pfadd@cke

\newcount\@tempcnta
\newcount\@tempcntb

\newdimen\@tempdima
\newdimen\@tempdimb

\newdimen\@wholewidth
\newdimen\@halfwidth

\newcount\@xarg
\newcount\@yarg
\newcount\@yyarg
\newbox\@linechar
\newbox\@tempboxa
\newdimen\@linelen
\newdimen\@clnwd
\newdimen\@clnht

\newif\if@negarg

\def\@whilenoop#1{}
\def\@whiledim#1\do #2{\ifdim #1\relax#2\@iwhiledim{#1\relax#2}\fi}
\def\@iwhiledim#1{\ifdim #1\let\@nextwhile=\@iwhiledim
        \else\let\@nextwhile=\@whilenoop\fi\@nextwhile{#1}}

\def\@whileswnoop#1\fi{}
\def\@whilesw#1\fi#2{#1#2\@iwhilesw{#1#2}\fi\fi}
\def\@iwhilesw#1\fi{#1\let\@nextwhile=\@iwhilesw
         \else\let\@nextwhile=\@whileswnoop\fi\@nextwhile{#1}\fi}

\def\thinlines{\let\@linefnt\tenln \let\@circlefnt\tencirc
  \@wholewidth\fontdimen8\tenln \@halfwidth .5\@wholewidth}
\def\thicklines{\let\@linefnt\tenlnw \let\@circlefnt\tencircw
  \@wholewidth\fontdimen8\tenlnw \@halfwidth .5\@wholewidth}
\thinlines

\PfadD@cke1pt \Pfadd@cke0.5pt
\def\PfadDicke#1{\PfadD@cke#1 \divide\PfadD@cke by2 \Pfadd@cke\PfadD@cke \multiply\PfadD@cke by2}
\long\def\LOOP#1\REPEAT{\def\BODY{#1}\ITERATE}
\def\ITERATE{\BODY \let\next\ITERATE \else\let\next\relax\fi \next}
\let\REPEAT=\fi
\def\Punkt{\hbox{\raise-2pt\hbox to0pt{\hss$\ssize\bullet$\hss}}}
\def\DuennPunkt(#1,#2){\unskip
  \raise#2 \Einheit\hbox to0pt{\hskip#1 \Einheit
          \raise-2.5pt\hbox to0pt{\hss$\bullet$\hss}\hss}}
\def\NormalPunkt(#1,#2){\unskip
  \raise#2 \Einheit\hbox to0pt{\hskip#1 \Einheit
          \raise-3pt\hbox to0pt{\hss\twelvepoint$\bullet$\hss}\hss}}
\def\DickPunkt(#1,#2){\unskip
  \raise#2 \Einheit\hbox to0pt{\hskip#1 \Einheit
          \raise-4pt\hbox to0pt{\hss\fourteenpoint$\bullet$\hss}\hss}}
\def\Kreis(#1,#2){\unskip
  \raise#2 \Einheit\hbox to0pt{\hskip#1 \Einheit
          \raise-4pt\hbox to0pt{\hss\fourteenpoint$\circ$\hss}\hss}}

\def\Line@(#1,#2)#3{\@xarg #1\relax \@yarg #2\relax
\@linelen=#3\Einheit
\ifnum\@xarg =0 \@vline
  \else \ifnum\@yarg =0 \@hline \else \@sline\fi
\fi}

\def\@sline{\ifnum\@xarg< 0 \@negargtrue \@xarg -\@xarg \@yyarg -\@yarg
  \else \@negargfalse \@yyarg \@yarg \fi
\ifnum \@yyarg >0 \@tempcnta\@yyarg \else \@tempcnta -\@yyarg \fi
\ifnum\@tempcnta>6 \@badlinearg\@tempcnta0 \fi
\ifnum\@xarg>6 \@badlinearg\@xarg 1 \fi
\setbox\@linechar\hbox{\@linefnt\@getlinechar(\@xarg,\@yyarg)}%
\ifnum \@yarg >0 \let\@upordown\raise \@clnht\z@
   \else\let\@upordown\lower \@clnht \ht\@linechar\fi
\@clnwd=\wd\@linechar
\if@negarg \hskip -\wd\@linechar \def\@tempa{\hskip -2\wd\@linechar}\else
     \let\@tempa\relax \fi
\@whiledim \@clnwd <\@linelen \do
  {\@upordown\@clnht\copy\@linechar
   \@tempa
   \advance\@clnht \ht\@linechar
   \advance\@clnwd \wd\@linechar}%
\advance\@clnht -\ht\@linechar
\advance\@clnwd -\wd\@linechar
\@tempdima\@linelen\advance\@tempdima -\@clnwd
\@tempdimb\@tempdima\advance\@tempdimb -\wd\@linechar
\if@negarg \hskip -\@tempdimb \else \hskip \@tempdimb \fi
\multiply\@tempdima \@m
\@tempcnta \@tempdima \@tempdima \wd\@linechar \divide\@tempcnta \@tempdima
\@tempdima \ht\@linechar \multiply\@tempdima \@tempcnta
\divide\@tempdima \@m
\advance\@clnht \@tempdima
\ifdim \@linelen <\wd\@linechar
   \hskip \wd\@linechar
  \else\@upordown\@clnht\copy\@linechar\fi}

\def\@hline{\ifnum \@xarg <0 \hskip -\@linelen \fi
\vrule height\Pfadd@cke width \@linelen depth\Pfadd@cke
\ifnum \@xarg <0 \hskip -\@linelen \fi}

\def\@getlinechar(#1,#2){\@tempcnta#1\relax\multiply\@tempcnta 8
\advance\@tempcnta -9 \ifnum #2>0 \advance\@tempcnta #2\relax\else
\advance\@tempcnta -#2\relax\advance\@tempcnta 64 \fi
\char\@tempcnta}

\def\Vektor(#1,#2)#3(#4,#5){\unskip\leavevmode
  \xcoord#4\relax \ycoord#5\relax
      \raise\ycoord \Einheit\hbox to0pt{\hskip\xcoord \Einheit
         \Vector@(#1,#2){#3}\hss}}

\def\Vector@(#1,#2)#3{\@xarg #1\relax \@yarg #2\relax
\@tempcnta \ifnum\@xarg<0 -\@xarg\else\@xarg\fi
\ifnum\@tempcnta<5\relax
\@linelen=#3\Einheit
\ifnum\@xarg =0 \@vvector
  \else \ifnum\@yarg =0 \@hvector \else \@svector\fi
\fi
\else\@badlinearg\fi}

\def\@hvector{\@hline\hbox to 0pt{\@linefnt
\ifnum \@xarg <0 \@getlarrow(1,0)\hss\else
    \hss\@getrarrow(1,0)\fi}}

\def\@vvector{\ifnum \@yarg <0 \@downvector \else \@upvector \fi}

\def\@svector{\@sline
\@tempcnta\@yarg \ifnum\@tempcnta <0 \@tempcnta=-\@tempcnta\fi
\ifnum\@tempcnta <5
  \hskip -\wd\@linechar
  \@upordown\@clnht \hbox{\@linefnt  \if@negarg
  \@getlarrow(\@xarg,\@yyarg) \else \@getrarrow(\@xarg,\@yyarg) \fi}%
\else\@badlinearg\fi}

\def\@upline{\hbox to \z@{\hskip -.5\Pfadd@cke \vrule width \Pfadd@cke
   height \@linelen depth \z@\hss}}

\def\@downline{\hbox to \z@{\hskip -.5\Pfadd@cke \vrule width \Pfadd@cke
   height \z@ depth \@linelen \hss}}

\def\@upvector{\@upline\setbox\@tempboxa\hbox{\@linefnt\char'66}\raise
     \@linelen \hbox to\z@{\lower \ht\@tempboxa\box\@tempboxa\hss}}

\def\@downvector{\@downline\lower \@linelen
      \hbox to \z@{\@linefnt\char'77\hss}}

\def\@getlarrow(#1,#2){\ifnum #2 =\z@ \@tempcnta='33\else
\@tempcnta=#1\relax\multiply\@tempcnta \sixt@@n \advance\@tempcnta
-9 \@tempcntb=#2\relax\multiply\@tempcntb \tw@
\ifnum \@tempcntb >0 \advance\@tempcnta \@tempcntb\relax
\else\advance\@tempcnta -\@tempcntb\advance\@tempcnta 64
\fi\fi\char\@tempcnta}

\def\@getrarrow(#1,#2){\@tempcntb=#2\relax
\ifnum\@tempcntb < 0 \@tempcntb=-\@tempcntb\relax\fi
\ifcase \@tempcntb\relax \@tempcnta='55 \or
\ifnum #1<3 \@tempcnta=#1\relax\multiply\@tempcnta
24 \advance\@tempcnta -6 \else \ifnum #1=3 \@tempcnta=49
\else\@tempcnta=58 \fi\fi\or
\ifnum #1<3 \@tempcnta=#1\relax\multiply\@tempcnta
24 \advance\@tempcnta -3 \else \@tempcnta=51\fi\or
\@tempcnta=#1\relax\multiply\@tempcnta
\sixt@@n \advance\@tempcnta -\tw@ \else
\@tempcnta=#1\relax\multiply\@tempcnta
\sixt@@n \advance\@tempcnta 7 \fi\ifnum #2<0 \advance\@tempcnta 64 \fi
\char\@tempcnta}

\def\Diagonale(#1,#2)#3{\unskip\leavevmode
  \xcoord#1\relax \ycoord#2\relax
      \raise\ycoord \Einheit\hbox to0pt{\hskip\xcoord \Einheit
         \Line@(1,1){#3}\hss}}
\def\AntiDiagonale(#1,#2)#3{\unskip\leavevmode
  \xcoord#1\relax \ycoord#2\relax 
      \raise\ycoord \Einheit\hbox to0pt{\hskip\xcoord \Einheit
         \Line@(1,-1){#3}\hss}}
\def\Pfad(#1,#2),#3\endPfad{\unskip\leavevmode
  \xcoord#1 \ycoord#2 \thicklines\ZeichnePfad#3\endPfad\thinlines}
\def\ZeichnePfad#1{\ifx#1\endPfad\let\next\relax
  \else\let\next\ZeichnePfad
    \ifnum#1=1
      \raise\ycoord \Einheit\hbox to0pt{\hskip\xcoord \Einheit
         \vrule height\Pfadd@cke width1 \Einheit depth\Pfadd@cke\hss}%
      \advance\xcoord by 1
    \else\ifnum#1=2
      \raise\ycoord \Einheit\hbox to0pt{\hskip\xcoord \Einheit
        \hbox{\hskip-\PfadD@cke\vrule height1 \Einheit width\PfadD@cke depth0pt}\hss}%
      \advance\ycoord by 1
    \else\ifnum#1=3
      \raise\ycoord \Einheit\hbox to0pt{\hskip\xcoord \Einheit
         \Line@(1,1){1}\hss}
      \advance\xcoord by 1
      \advance\ycoord by 1
    \else\ifnum#1=4
      \raise\ycoord \Einheit\hbox to0pt{\hskip\xcoord \Einheit
         \Line@(1,-1){1}\hss}
      \advance\xcoord by 1
      \advance\ycoord by -1
    \fi\fi\fi\fi
  \fi\next}
\def\hSSchritt{\leavevmode\raise-.4pt\hbox to0pt{\hss.\hss}\hskip.2\Einheit
  \raise-.4pt\hbox to0pt{\hss.\hss}\hskip.2\Einheit
  \raise-.4pt\hbox to0pt{\hss.\hss}\hskip.2\Einheit
  \raise-.4pt\hbox to0pt{\hss.\hss}\hskip.2\Einheit
  \raise-.4pt\hbox to0pt{\hss.\hss}\hskip.2\Einheit}
\def\vSSchritt{\vbox{\baselineskip.2\Einheit\lineskiplimit0pt
\hbox{.}\hbox{.}\hbox{.}\hbox{.}\hbox{.}}}
\def\DSSchritt{\leavevmode\raise-.4pt\hbox to0pt{%
  \hbox to0pt{\hss.\hss}\hskip.2\Einheit
  \raise.2\Einheit\hbox to0pt{\hss.\hss}\hskip.2\Einheit
  \raise.4\Einheit\hbox to0pt{\hss.\hss}\hskip.2\Einheit
  \raise.6\Einheit\hbox to0pt{\hss.\hss}\hskip.2\Einheit
  \raise.8\Einheit\hbox to0pt{\hss.\hss}\hss}}
\def\dSSchritt{\leavevmode\raise-.4pt\hbox to0pt{%
  \hbox to0pt{\hss.\hss}\hskip.2\Einheit
  \raise-.2\Einheit\hbox to0pt{\hss.\hss}\hskip.2\Einheit
  \raise-.4\Einheit\hbox to0pt{\hss.\hss}\hskip.2\Einheit
  \raise-.6\Einheit\hbox to0pt{\hss.\hss}\hskip.2\Einheit
  \raise-.8\Einheit\hbox to0pt{\hss.\hss}\hss}}
\def\SPfad(#1,#2),#3\endSPfad{\unskip\leavevmode
  \xcoord#1 \ycoord#2 \ZeichneSPfad#3\endSPfad}
\def\ZeichneSPfad#1{\ifx#1\endSPfad\let\next\relax
  \else\let\next\ZeichneSPfad
    \ifnum#1=1
      \raise\ycoord \Einheit\hbox to0pt{\hskip\xcoord \Einheit
         \hSSchritt\hss}%
      \advance\xcoord by 1
    \else\ifnum#1=2
      \raise\ycoord \Einheit\hbox to0pt{\hskip\xcoord \Einheit
        \hbox{\hskip-2pt \vSSchritt}\hss}%
      \advance\ycoord by 1
    \else\ifnum#1=3
      \raise\ycoord \Einheit\hbox to0pt{\hskip\xcoord \Einheit
         \DSSchritt\hss}
      \advance\xcoord by 1
      \advance\ycoord by 1
    \else\ifnum#1=4
      \raise\ycoord \Einheit\hbox to0pt{\hskip\xcoord \Einheit
         \dSSchritt\hss}
      \advance\xcoord by 1
      \advance\ycoord by -1
    \fi\fi\fi\fi
  \fi\next}
\def\Koordinatenachsen(#1,#2){\unskip
 \hbox to0pt{\hskip-.5pt\vrule height#2 \Einheit width.5pt depth1 \Einheit}%
 \hbox to0pt{\hskip-1 \Einheit \xcoord#1 \advance\xcoord by1
    \vrule height0.25pt width\xcoord \Einheit depth0.25pt\hss}}
\def\Koordinatenachsen(#1,#2)(#3,#4){\unskip
 \hbox to0pt{\hskip-.5pt \ycoord-#4 \advance\ycoord by1
    \vrule height#2 \Einheit width.5pt depth\ycoord \Einheit}%
 \hbox to0pt{\hskip-1 \Einheit \hskip#3\Einheit 
    \xcoord#1 \advance\xcoord by1 \advance\xcoord by-#3 
    \vrule height0.25pt width\xcoord \Einheit depth0.25pt\hss}}
\def\Gitter(#1,#2){\unskip \xcoord0 \ycoord0 \leavevmode
  \LOOP\ifnum\ycoord<#2
    \loop\ifnum\xcoord<#1
      \raise\ycoord \Einheit\hbox to0pt{\hskip\xcoord \Einheit\Punkt\hss}%
      \advance\xcoord by1
    \repeat
    \xcoord0
    \advance\ycoord by1
  \REPEAT}
\def\Gitter(#1,#2)(#3,#4){\unskip \xcoord#3 \ycoord#4 \leavevmode
  \LOOP\ifnum\ycoord<#2
    \loop\ifnum\xcoord<#1
      \raise\ycoord \Einheit\hbox to0pt{\hskip\xcoord \Einheit\Punkt\hss}%
      \advance\xcoord by1
    \repeat
    \xcoord#3
    \advance\ycoord by1
  \REPEAT}
\def\Label#1#2(#3,#4){\unskip \xdim#3 \Einheit \ydim#4 \Einheit
  \def\lo{\advance\xdim by-.5 \Einheit \advance\ydim by.5 \Einheit}%
  \def\llo{\advance\xdim by-.25cm \advance\ydim by.5 \Einheit}%
  \def\loo{\advance\xdim by-.5 \Einheit \advance\ydim by.25cm}%
  \def\o{\advance\ydim by.25cm}%
  \def\ro{\advance\xdim by.5 \Einheit \advance\ydim by.5 \Einheit}%
  \def\rro{\advance\xdim by.25cm \advance\ydim by.5 \Einheit}%
  \def\roo{\advance\xdim by.5 \Einheit \advance\ydim by.25cm}%
  \def\l{\advance\xdim by-.30cm}%
  \def\r{\advance\xdim by.30cm}%
  \def\lu{\advance\xdim by-.5 \Einheit \advance\ydim by-.6 \Einheit}%
  \def\llu{\advance\xdim by-.25cm \advance\ydim by-.6 \Einheit}%
  \def\luu{\advance\xdim by-.5 \Einheit \advance\ydim by-.30cm}%
  \def\u{\advance\ydim by-.30cm}%
  \def\ru{\advance\xdim by.5 \Einheit \advance\ydim by-.6 \Einheit}%
  \def\rru{\advance\xdim by.25cm \advance\ydim by-.6 \Einheit}%
  \def\ruu{\advance\xdim by.5 \Einheit \advance\ydim by-.30cm}%
  #1\raise\ydim\hbox to0pt{\hskip\xdim
     \vbox to0pt{\vss\hbox to0pt{\hss$#2$\hss}\vss}\hss}%
}
\catcode`\@=13

\catcode`\@=11
\font@\twelverm=cmr10 scaled\magstep1
\font@\twelveit=cmti10 scaled\magstep1
\font@\twelvebf=cmbx10 scaled\magstep1
\font@\twelvei=cmmi10 scaled\magstep1
\font@\twelvesy=cmsy10 scaled\magstep1
\font@\twelveex=cmex10 scaled\magstep1

\newtoks\twelvepoint@
\def\twelvepoint{\normalbaselineskip15\p@
 \abovedisplayskip15\p@ plus3.6\p@ minus10.8\p@
 \belowdisplayskip\abovedisplayskip
 \abovedisplayshortskip\z@ plus3.6\p@
 \belowdisplayshortskip8.4\p@ plus3.6\p@ minus4.8\p@
 \textonlyfont@\rm\twelverm \textonlyfont@\it\twelveit
 \textonlyfont@\sl\twelvesl \textonlyfont@\bf\twelvebf
 \textonlyfont@\smc\twelvesmc \textonlyfont@\tt\twelvett
%
 \ifsyntax@ \def\big##1{{\hbox{$\left##1\right.$}}}%
  \let\Big\big \let\bigg\big \let\Bigg\big
 \else
  \textfont\z@=\twelverm  \scriptfont\z@=\tenrm  \scriptscriptfont\z@=\sevenrm
  \textfont\@ne=\twelvei  \scriptfont\@ne=\teni  \scriptscriptfont\@ne=\seveni
  \textfont\tw@=\twelvesy \scriptfont\tw@=\tensy \scriptscriptfont\tw@=\sevensy
  \textfont\thr@@=\twelveex \scriptfont\thr@@=\tenex
        \scriptscriptfont\thr@@=\tenex
  \textfont\itfam=\twelveit \scriptfont\itfam=\tenit
        \scriptscriptfont\itfam=\tenit
  \textfont\bffam=\twelvebf \scriptfont\bffam=\tenbf
        \scriptscriptfont\bffam=\sevenbf
  \setbox\strutbox\hbox{\vrule height10.2\p@ depth4.2\p@ width\z@}%
  \setbox\strutbox@\hbox{\lower.6\normallineskiplimit\vbox{%
        \kern-\normallineskiplimit\copy\strutbox}}%
 \setbox\z@\vbox{\hbox{$($}\kern\z@}\bigsize@=1.4\ht\z@
 \fi
 \normalbaselines\rm\ex@.2326ex\jot3.6\ex@\the\twelvepoint@}

\font@\fourteenrm=cmr10 scaled\magstep2
\font@\fourteenit=cmti10 scaled\magstep2
\font@\fourteensl=cmsl10 scaled\magstep2
\font@\fourteensmc=cmcsc10 scaled\magstep2
\font@\fourteentt=cmtt10 scaled\magstep2
\font@\fourteenbf=cmbx10 scaled\magstep2
\font@\fourteeni=cmmi10 scaled\magstep2
\font@\fourteensy=cmsy10 scaled\magstep2
\font@\fourteenex=cmex10 scaled\magstep2
\font@\fourteenmsa=msam10 scaled\magstep2
\font@\fourteeneufm=eufm10 scaled\magstep2
\font@\fourteenmsb=msbm10 scaled\magstep2
\newtoks\fourteenpoint@
\def\fourteenpoint{\normalbaselineskip15\p@
 \abovedisplayskip18\p@ plus4.3\p@ minus12.9\p@
 \belowdisplayskip\abovedisplayskip
 \abovedisplayshortskip\z@ plus4.3\p@
 \belowdisplayshortskip10.1\p@ plus4.3\p@ minus5.8\p@
 \textonlyfont@\rm\fourteenrm \textonlyfont@\it\fourteenit
 \textonlyfont@\sl\fourteensl \textonlyfont@\bf\fourteenbf
 \textonlyfont@\smc\fourteensmc \textonlyfont@\tt\fourteentt
%
 \ifsyntax@ \def\big##1{{\hbox{$\left##1\right.$}}}%
  \let\Big\big \let\bigg\big \let\Bigg\big
 \else
  \textfont\z@=\fourteenrm  \scriptfont\z@=\twelverm  \scriptscriptfont\z@=\tenrm
  \textfont\@ne=\fourteeni  \scriptfont\@ne=\twelvei  \scriptscriptfont\@ne=\teni
  \textfont\tw@=\fourteensy \scriptfont\tw@=\twelvesy \scriptscriptfont\tw@=\tensy
  \textfont\thr@@=\fourteenex \scriptfont\thr@@=\twelveex
        \scriptscriptfont\thr@@=\twelveex
  \textfont\itfam=\fourteenit \scriptfont\itfam=\twelveit
        \scriptscriptfont\itfam=\twelveit
  \textfont\bffam=\fourteenbf \scriptfont\bffam=\twelvebf
        \scriptscriptfont\bffam=\tenbf
  \setbox\strutbox\hbox{\vrule height12.2\p@ depth5\p@ width\z@}%
  \setbox\strutbox@\hbox{\lower.72\normallineskiplimit\vbox{%
        \kern-\normallineskiplimit\copy\strutbox}}%
 \setbox\z@\vbox{\hbox{$($}\kern\z@}\bigsize@=1.7\ht\z@
 \fi
 \normalbaselines\rm\ex@.2326ex\jot4.3\ex@\the\fourteenpoint@}

\font@\seventeenrm=cmr10 scaled\magstep3
\font@\seventeenit=cmti10 scaled\magstep3
\font@\seventeensl=cmsl10 scaled\magstep3
\font@\seventeensmc=cmcsc10 scaled\magstep3
\font@\seventeentt=cmtt10 scaled\magstep3
\font@\seventeenbf=cmbx10 scaled\magstep3
\font@\seventeeni=cmmi10 scaled\magstep3
\font@\seventeensy=cmsy10 scaled\magstep3
\font@\seventeenex=cmex10 scaled\magstep3
\font@\seventeenmsa=msam10 scaled\magstep3
\font@\seventeeneufm=eufm10 scaled\magstep3
\font@\seventeenmsb=msbm10 scaled\magstep3
\newtoks\seventeenpoint@
\def\seventeenpoint{\normalbaselineskip18\p@
 \abovedisplayskip21.6\p@ plus5.2\p@ minus15.4\p@
 \belowdisplayskip\abovedisplayskip
 \abovedisplayshortskip\z@ plus5.2\p@
 \belowdisplayshortskip12.1\p@ plus5.2\p@ minus7\p@
 \textonlyfont@\rm\seventeenrm \textonlyfont@\it\seventeenit
 \textonlyfont@\sl\seventeensl \textonlyfont@\bf\seventeenbf
 \textonlyfont@\smc\seventeensmc \textonlyfont@\tt\seventeentt
%
 \ifsyntax@ \def\big##1{{\hbox{$\left##1\right.$}}}%
  \let\Big\big \let\bigg\big \let\Bigg\big
 \else
  \textfont\z@=\seventeenrm  \scriptfont\z@=\fourteenrm  \scriptscriptfont\z@=\twelverm
  \textfont\@ne=\seventeeni  \scriptfont\@ne=\fourteeni  \scriptscriptfont\@ne=\twelvei
  \textfont\tw@=\seventeensy \scriptfont\tw@=\fourteensy \scriptscriptfont\tw@=\twelvesy
  \textfont\thr@@=\seventeenex \scriptfont\thr@@=\fourteenex
        \scriptscriptfont\thr@@=\fourteenex
  \textfont\itfam=\seventeenit \scriptfont\itfam=\fourteenit
        \scriptscriptfont\itfam=\fourteenit
  \textfont\bffam=\seventeenbf \scriptfont\bffam=\fourteenbf
        \scriptscriptfont\bffam=\twelvebf
  \setbox\strutbox\hbox{\vrule height14.6\p@ depth6\p@ width\z@}%
  \setbox\strutbox@\hbox{\lower.86\normallineskiplimit\vbox{%
        \kern-\normallineskiplimit\copy\strutbox}}%
 \setbox\z@\vbox{\hbox{$($}\kern\z@}\bigsize@=2\ht\z@
 \fi
 \normalbaselines\rm\ex@.2326ex\jot5.2\ex@\the\seventeenpoint@}

\catcode`\@=13

\def\ldreieck{\bsegment
  \rlvec(1 0) \rlvec(-.5 -0.866025403784439)
  \rlvec(-.5 0.866025403784439)  
  \savepos(.5 -0.866025403784439)(*ex *ey)
        \esegment
  \move(*ex *ey)
        }
\def\rdreieck{\bsegment
  \rlvec(.5 -0.866025403784439) \rlvec(-1 0)  
  \rlvec(.5 0.866025403784439)
  \savepos(-.5 -0.866025403784439)(*ex *ey)
        \esegment
  \move(*ex *ey)
        }
\def\rhombus{\bsegment
  \rlvec(1 0) \rlvec(.5 -0.866025403784439) 
  \rlvec(-1 0)  \rlvec(.5 0.866025403784439)  
  \rmove(-.5 -0.866025403784439)  \rlvec(-.5 0.866025403784439) 
  \savepos(.5 -0.866025403784439)(*ex *ey)
        \esegment
  \move(*ex *ey)
        }
\def\RhombusA{\bsegment
  \rlvec(1 0) \rlvec(.5 -0.866025403784439) 
  \rlvec(-1 0) \rlvec(-.5 0.866025403784439) 
  \savepos(.5 -0.866025403784439)(*ex *ey)
        \esegment
  \move(*ex *ey)
        }
\def\RhombusB{\bsegment
  \rlvec(1 0) \rlvec(-.5 -0.866025403784439)
  \rlvec(-1 0) \rlvec(.5 0.866025403784439) 
  \savepos(-.5 -0.866025403784439)(*ex *ey)
        \esegment
  \move(*ex *ey)
        }
\def\RhombusC{\bsegment
  \rlvec(.5 -0.866025403784439) \rlvec(-.5 -0.866025403784439)
  \rlvec(-.5 0.866025403784439) \rlvec(.5 0.866025403784439) 
  \savepos(.5 -0.866025403784439)(*ex *ey)
        \esegment
  \move(*ex *ey)
        }

\def\hdSchritt{\bsegment
  \lpatt(.05 .13)
  \rlvec(.5 -0.866025403784439) 
  \savepos(.5 -0.866025403784439)(*ex *ey)
        \esegment
  \move(*ex *ey)
        }

\def\vdaSchritt{\bsegment
  \lpatt(.05 .13)
  \rlvec(.5 0.866025403784439) 
  \savepos(.5 0.866025403784439)(*ex *ey)
        \esegment
  \move(*ex *ey)
        }

\def\odaSchritt{\bsegment
  \lpatt(.05 .13)
  \rlvec(1 0) 
  \savepos(1 0)(*ex *ey)
        \esegment
  \move(*ex *ey)
        }

\def\ringerl(#1 #2){\move(#1 #2)\fcir f:0 r:.1}

\topmatter 
\title A non-automatic (!) application of Gosper's algorithm evaluates
a determinant from tiling enumeration
\endtitle 
\author M.~Ciucu and C.~Krattenthaler\footnote"$^\dagger$"{\hbox{Partially supported by the Austrian
Science Foundation FWF, grant P13190-MAT.}}
\endauthor 
\affil Department of Mathematics, \\
Georgia Institute of Technology,\\
Atlanta, GA 30332-0160, USA\\
e-mail: {\tt ciucu\@math.gatech.edu}\\
WWW: {\tt http://www.math.gatech.edu/\~{}ciucu}\\\vskip6pt
Institut f\"ur Mathematik der Universit\"at Wien,\\
Strudlhofgasse 4, A-1090 Wien, Austria.\\
e-mail: KRATT\@Ap.Univie.Ac.At\\
WWW: \tt http://www.mat.univie.ac.at/People/kratt
\endaffil
\address Department of Mathematics, 
Georgia Institute of Technology,
Atlanta, GA 30332-0160, USA.
\endaddress
\address Institut f\"ur Mathematik der Universit\"at Wien,
Strudlhofgasse 4, A-1090 Wien, Austria.
\endaddress
\abstract We evaluate the determinant 
$\det_{1\leq
i,j\leq n}\left(\binom{x+y+j}{x-i+2j}-\binom{x+y+j}{x+i+2j}\right)$, 
which gives the
number of lozenge tilings of a hexagon with cut off corners. A
particularly interesting feature of this evaluation is that it
requires the proof of a certain hypergeometric identity which we
accomplish by using Gosper's algorithm in a non-automatic fashion.
\endabstract
\endtopmatter

\rightheadtext{Gosper's algorithm evaluates a determinant}

\leftheadtext{M. Ciucu and C. Krattenthaler}

\document

The purpose of this paper is to provide a direct evaluation of the
determinant
$$\det_{1\le i,j\le n}\(\binom {x+y+j}{x-i+2j}-\binom
{x+y+j}{x+i+2j}\).\tag\ZAa$$
This determinant arises in our study \cite{\CiKrAD} on
the enumeration of lozenge tilings of hexagons with cut off corners.
For example, consider a hexagon with side lengths $x+n$, $n$, $y$,
$x+n$, $n$, $y$ (in cyclic order) and angles of $120^\circ$ of which
two adjacent corners are cut off as in Figure~1(a).\footnote{To be precise,
from the top-left corner we cut off a (reversed) staircase of the form
$(y-1,y-2,\dots,1)$, meaning that the cut-off staircase consists of $y-1$
rhombi in the first row, $y-2$ rhombi in the second row, etc., and
from the top-right corner we cut off a staircase of the form
$(n-1,n-2,\dots,1)$.} Figure~1(b) shows a
lozenge tiling of this region, by which we mean 
a tiling by unit rhombi with angles of
$60^\circ$ and $120^\circ$, referred to as lozenges. The number of
these
lozenge tilings is given by the determinant (\ZAa). This is seen by
converting the lozenge tilings into families $(P_1,P_2,\dots,P_n)$ of
nonintersecting lattice paths consisting of positive unit steps, where
the path $P_i$ runs from $(i,-i)$ to $(x+2i,y-i)$, $i=1,2,\dots,n$
and does not cross the diagonal $y=x-1$ (see Figure~2), 
and then applying the
main theorem of nonintersecting lattice paths 
\cite{\LindAA, Lemma~1}, \cite{\GeViAA}, \cite{\StemAE, Theorem~1.2} (see
\cite{\CiKrAD} for details and background; there is also another case in
\cite{\CiKrAD} in which the determinant (\ZAa) provides the solution).

\midinsert
\vskip15pt
\vbox{
\centertexdraw{
\drawdim truecm \linewd.02
\rhombus \rhombus \rhombus \rhombus \rhombus \ldreieck
\move(1.5 -.866025)
\rhombus \rhombus \rhombus \ldreieck
\move(3 -1.7305)
\rhombus \ldreieck
\move(-1 0)
\rhombus \rhombus \rhombus \rhombus \rhombus \rhombus \ldreieck
\move(-1 0)
\rdreieck \rhombus \rhombus \rhombus \rhombus \rhombus \rhombus
\move(-2.5 -.866025)
\rhombus \rhombus \rhombus \rhombus \rhombus \rhombus 
\move(-2.5 -.866025)
\rdreieck \rhombus \rhombus \rhombus \rhombus \rhombus 
\move(-4 -1.7305)
\rhombus \rhombus \rhombus \rhombus \rhombus 
\move(-4 -1.7305)
\rdreieck \rhombus \rhombus \rhombus \rhombus 
\move(-5.5 -2.596)
\rhombus \rhombus \rhombus \rhombus 
\move(-5.5 -2.596)
\rdreieck \rhombus \rhombus \rhombus 
\move (-5.5 -2.596)
\vdaSchritt \vdaSchritt \vdaSchritt 
\odaSchritt \odaSchritt \odaSchritt 
\move (1 0)
\odaSchritt \odaSchritt \hdSchritt \hdSchritt 

\htext(-1 0.4){$x+n$}
\htext(-1.5 -6.7){$x+n$}
\htext(-5.6 -1.8){$y$}
\htext(3.9 -4.6){$y$}
\htext(4.2 -1.3){$n$}
\htext(-6 -5.1){$n$}
\rtext td:30 (-5.6 -5){$\left\{ \vbox{\vskip1.4cm} \right.$}
\rtext td:30 (3.8 -1.4){$\left. \vbox{\vskip1.4cm} \right\}$}
\rtext td:90 (-.4 -.1) {$\left. \vbox{\vskip3.2cm} \right\} $}
\rtext td:90 (-0.9 -6.4){$\left\{ \vbox{\vskip3.2cm}\right. $}
\rtext td:-30 (-5.4 -1.7){$\left\{ \vbox{\vskip1.8cm} \right. $}
\rtext td:-30 (3.4 -4.4) {$\left. \vbox{\vskip1.8cm} \right\} $}
}
\vskip4pt
\centerline{\eightpoint (a) A hexagon with cut off corners.}
\vskip10pt

\centertexdraw{
\drawdim truecm \linewd.05
\move (12 0)
\RhombusC \RhombusB 
\rmove (1.5 .86602)
\RhombusC \RhombusC \RhombusB \RhombusB \RhombusB \RhombusB
\move (14 -1.7305)
\RhombusA
\move (11 0)
\RhombusB \RhombusA \RhombusA \RhombusB \RhombusA \RhombusB \RhombusB
\move (10 0)
\RhombusB \RhombusA \RhombusA \RhombusB \RhombusB \RhombusA \RhombusB
\move (12 -1.7305)
\RhombusA \RhombusB \RhombusA \RhombusB \RhombusB 
\move (8.5 -.866025)
\RhombusA \RhombusB \RhombusA \RhombusB \RhombusA \RhombusB
\move (8.5 -.866025)
\RhombusC 
\rmove (-1 0)
\RhombusC
\move (7 -1.7305)
\RhombusB \RhombusA \RhombusA \RhombusA \RhombusB
\move (5.5 -2.596)
\RhombusA \RhombusB \RhombusA \RhombusA 
\move (7 -3.4641)
\RhombusC
\move (5.5 -2.596)
\RhombusC
}
\vskip4pt
\centerline{\eightpoint (b) A lozenge tiling of the 
hexagon with cut off corners.}
\vskip10pt
\centerline{\eightpoint Figure 1}
}
\vskip10pt
\endinsert

\midinsert
\vskip15pt
\vbox{
\centertexdraw{
\drawdim truecm \linewd.05
\move (12 0)
\RhombusC \RhombusB 
\rmove (1.5 .86602)
\RhombusC \RhombusC \RhombusB \RhombusB \RhombusB \RhombusB
\move (14 -1.7305)
\RhombusA
\move (11 0)
\RhombusB \RhombusA \RhombusA \RhombusB \RhombusA \RhombusB \RhombusB
\move (10 0)
\RhombusB \RhombusA \RhombusA \RhombusB \RhombusB \RhombusA \RhombusB
\move (12 -1.7305)
\RhombusA \RhombusB \RhombusA \RhombusB \RhombusB 
\move (8.5 -.866025)
\RhombusA \RhombusB \RhombusA \RhombusB \RhombusA \RhombusB
\move (8.5 -.866025)
\RhombusC 
\rmove (-1 0)
\RhombusC
\move (7 -1.7305)
\RhombusB \RhombusA \RhombusA \RhombusA \RhombusB
\move (5.5 -2.596)
\RhombusA \RhombusB \RhombusA \RhombusA 
\move (7 -3.4641)
\RhombusC
\move (5.5 -2.596)
\RhombusC

\ringerl (5.25 -3.89711)
\vdaSchritt \odaSchritt \odaSchritt \vdaSchritt \vdaSchritt 
\odaSchritt \odaSchritt \odaSchritt \vdaSchritt 
\ringerl (12.25 -.433012)
\ringerl (5.75 -4.76314)
\odaSchritt \vdaSchritt \odaSchritt \vdaSchritt \odaSchritt
\vdaSchritt \odaSchritt \odaSchritt \odaSchritt \vdaSchritt  
\ringerl (13.75 -1.29904)
\ringerl (6.25 -5.62917)
\odaSchritt \vdaSchritt \odaSchritt \odaSchritt \odaSchritt
\vdaSchritt \odaSchritt \odaSchritt \vdaSchritt \vdaSchritt
\odaSchritt  
\ringerl (15.25 -2.16506)
}
\vskip4pt
\centerline{\eightpoint (a) The path family corresponding to a lozenge
tiling.}

$$
\Einheit=.7cm
\Gitter(11,4)(0,-3)
\Koordinatenachsen(11,4)(0,-3)
\Pfad(1,-1),211221112\endPfad
\Pfad(2,-2),1212121112\endPfad
\Pfad(3,-3),12111211212\endPfad
\Diagonale(-1,-2){6}
\DickPunkt(1,-1)
\DickPunkt(2,-2)
\DickPunkt(3,-3)
\DickPunkt(6,3)
\DickPunkt(8,2)
\DickPunkt(10,1)
\hskip7cm
$$
\vskip4pt
\centerline{\eightpoint (b) The paths made orthogonal.}
\vskip10pt

\vskip10pt
\centerline{\eightpoint Figure 2}
}
\vskip10pt
\endinsert

By Theorem~\TA\ below, the number of the lozenge tilings of the
preceding paragraph is given by a closed form expression.
The proof of Theorem~\TA\ that we present in this paper%
\footnote{An alternative proof is presented in \cite{\CiKrAD}, 
in which a combinatorial argument is used to convert the determinant (\ZAa) 
into a different determinant that was already known from
\cite{\KratBD, Theorem~10}.} is primarily based on
hypergeometric series identities. A remarkable aspect is that it
contains an instance of a {\it non-automatic}
application of Gosper's algorithm \cite{\GospAB} 
(see also 
\cite{\GrKPAA, \S5.7}, \cite{\PeWZAA, \S II.5}), see Step~3 of the
proof of Theorem~\TA. This is noteworthy,
because Gosper invented his algorithm to {\it automate} summation, 
so that a non-automatic application must be almost
considered as a misuse. But clearly (and more seriously), 
the fact that Gosper's algorithm
is also useful in ``computer-free territory'' only adds to its value.
(The only other instance of a non-automatic application of Gosper's
algorithm that we are aware of appears in \cite{\PeWiAB}. However,
the purpose of use there is different. Roughly speaking, we use it to
prove a {\it positive} result, namely to verify the truth of an
identity between certain hypergeometric series, see (14). 
In contrast, Petkov\v sek and Wilf use it to
prove a {\it negative} result, namely that a certain binomial sum cannot be
expressed in terms of closed form expressions.)

\proclaim{Theorem~\TA}Let $n$ be a positive integer, and let $x$ and
$y$ be nonnegative integers. Then the
following determinant evaluation holds:
$$\multline \det_{1\le i,j\le n}\(\binom {x+y+j}{x-i+2j}-\binom
{x+y+j}{x+i+2j}\)\\
=\prod _{j=1} ^{n}\frac {(j-1)!\,(x+y+2j)!\,(x-y+2j+1)_j\,
(x+2y+3j+1)_{n-j}} {(x+n+2j)!\,(y+n-j)!},
\endmultline\tag\ZA$$
where the shifted factorial
$(a)_k$ is defined by $(a)_k:=a(a+1)\cdots(a+k-1)$,
$k\ge1$, and $(a)_0:=1$.
\endproclaim

\remark{Remark} We formulate Theorem~\TA\ only for integral $x$ and
$y$. But in fact, with a
generalized definition of factorials and binomials (cf\.
\cite{\GrKPAA, \S5.5, (5.96), (5.100)},
Theorem~\TA\ would also make sense and be true for complex $x$ and $y$. 
\endremark

\demo{Proof} We prove the determinant evaluation by ``identification
of factors," a method that is also applied 
successfully in 
\cite{\CiEKAA}, \cite{\CiKrAA}, \cite{\EisTAC}, 
\cite{\FuKrAC}, \cite{\FuKrAD}, \cite{\KratBG}, \cite{\KratBD},
\cite{\KratBI},  
\cite{\KratBH}, \cite{\KrZeAA}, \cite{\KupeAD} and \cite{\PoorAB}
(see in particular the tutorial description in \cite{\KratBN, \S2.4} or
\cite{\KratBI, \S2}).

First of all, we take appropriate
factors out of the determinant. To be precise, we take 
$ {(x+y+j)!} /\big({(x+n+2j)!\,(y+n-j)!}\big)$ out of the
$j$-th column of the determinant in (\ZA), $j=1,2,\dots,n$. Thus we
obtain
$$\multline 
\prod _{i=1} ^{n}\frac {(x+y+j)!} {(x+n+2j)!\,(y+n-j)!}\\
\times \det_{1\le i,j\le n}\((x+2j-i+1)_{n+i}\,(y+i-j+1)_{n-i}-
(x+2j+i+1)_{n-i}\,(y-i-j+1)_{n+i}\)
\endmultline\tag\ZB$$
for the determinant in (\ZA).
Let us denote the determinant in (\ZB) by $D_n(x,y)$. 
Comparison of (\ZA) and (\ZB) yields that (\ZA) will be proved once
we are able to establish the determinant evaluation
$$\multline D_n(x,y)\\
=\det_{1\le i,j\le n}\((x+2j-i+1)_{n+i}\,(y+i-j+1)_{n-i}-
(x+2j+i+1)_{n-i}\,(y-i-j+1)_{n+i}\)\\
=\prod _{j=1} ^{n}{(j-1)!\,(x+y+j+1)_j\,(x-y+2j+1)_j\,
(x+2y+3j+1)_{n-j}} .
\endmultline\tag\ZC$$

For the proof of (\ZC) we proceed in several steps. An outline is as
follows. In the first step we show that $\prod _{j=1} ^{n}(x-y+2j+1)_j$ 
is a factor of $D_n(x,y)$  
as a polynomial in $x$ and $y$. In the second step we
show that $\prod _{j=1} ^{n}(x+y+j+1)_j$ 
is a factor of $D_n(x,y)$, and in the third step we
show that $\prod _{j=1} ^{n}(x+2y+3j+1)_{n-j}$ 
is a factor of $D_n(x,y)$. Then,
in the fourth step we determine the maximal degree of $D_n(x,y)$ as a
polynomial in $x$, and the maximal degree as a
polynomial in $y$, which turns out to be $n(3n+1)/2$ in both cases.
On the other hand, the degree in $x$, and also in $y$, of the product on
the right hand side of (\ZC), which by the first three steps divides
$D_n(x,y)$, is exactly $n(3n+1)/2$. Therefore we are forced to conclude that
$$D_n(x,y)=C(n)\prod _{j=1} ^{n}{(x-y+2j+1)_j\,(x+y+j+1)_j\,
(x+2y+3j+1)_{n-j}} ,\tag\ZD$$
where $C(n)$ is a constant independent of $x$ and $y$.
Finally, in the fifth step, we determine the constant $C(n)$, which
turns out to equal $\prod _{j=1} ^{n}(j-1)!$.
Clearly, this would finish the proof of (\ZC), and thus of (\ZA), 
as we already noted.

\smallskip
{\it Step 1. $\prod _{j=1} ^{n}(x-y+2j+1)_j$ is a factor of
$D_n(x,y)$}. Let us concentrate on a typical factor $(x-y+2j+l)$,
$1\le j\le n$, $1\le l\le j$. We claim that for each
such factor there is a linear combination of the columns that vanishes
if the factor vanishes. More precisely, we claim that for any $j,l$
with $1\le j\le n$, $1\le l\le j$ there holds
$$\multline 
 \sum_{s = l}^{\fl{{{j + l}\over 2}}}
\frac {\left( j - l \right)  } {\left( j - s \right) }
\frac {         ({ \textstyle j + l - 2 s +1}) _{s-l} } {\left( s-l \right) !}
     {{         ({ \textstyle x + 2 j + l + n - s + 1}) _{s-l} \,
         ({ \textstyle x + n + 2 s + 1}) _{j + l - 2 s} }\over 
       {         ({ \textstyle 2x + 2 j + l + s + 1}) _{j - s} }} \\
\cdot(\text {column
$s$ of $D_n(x,x+2j+l)$}) \\
+(\text {column
$j$ of $D_n(x,x+2j+l)$}) =0.
\endmultline\tag\ZE$$
To avoid confusion, for $j=l$ it is understood by convention that
the sum in (\ZE) vanishes.

In order to verify (\ZE), we have to check
$$\multline 
 \sum_{s = l}^{\fl{{{j + l}\over 2}}}
\frac {\left( j - l \right) } {\left( j - s \right) }
\frac {({ \textstyle j + l - 2 s+1}) _{s-l} } {\left( s-l \right) !}
     {{         ({ \textstyle x + 2 j + l + n - s + 1}) _{s-l} \,
         ({ \textstyle x + n + 2 s + 1}) _{j + l - 2 s} }\over 
       {         ({ \textstyle 2 x + 2 j + l + s + 1}) _{j - s} }} \\
\cdot         \big( ({ \textstyle x + i + 2 j + l - s + 1}) _{n-i} \,
            ({ \textstyle x - i + 2 s + 1}) _{n+i} \hskip5cm\\
\hskip3cm -       ({ \textstyle x - i + 2 j + l - s + 1}) _{n+i} \,
            ({ \textstyle x + i + 2 s + 1}) _{n-i}  \big)\\
+  ({ \textstyle x - i + 2 j + 1}) _{n+i} \,
   ({ \textstyle x + i + j + l + 1}) _{n-i} \hskip4cm\\
- ({ \textstyle x + i + 2 j + 1}) _{n-i} \,
     ({ \textstyle x - i + j + l + 1}) _{n+i}   =0,
\endmultline\tag\ZEa$$
which is (\ZE) restricted to the $i$-th row.
The exceptional case $j=l$ can be treated immediately. By
assumption, the sum in (\ZEa) vanishes for $j=l$, and, by
inspection, also 
the other two expressions in (\ZEa) vanish for $j=l$. So 
it remains to establish (\ZEa) for $j>l$.
In terms of the standard hypergeometric notation
$${}_r F_s\!\left[\matrix a_1,\dots,a_r\\ b_1,\dots,b_s\endmatrix; 
z\right]=\sum _{k=0} ^{\infty}\frac {\po{a_1}{k}\cdots\po{a_r}{k}}
{k!\,\po{b_1}{k}\cdots\po{b_s}{k}} z^k\ ,$$
this means to check
$$\multline 
  {{  
      ({ \textstyle x + i + 2 j + 1}) _{n-i} \,
      ({ \textstyle x - i + 2 l + 1}) _{n+i + j - l } }\over 
    {({ \textstyle 2 x + 2 j + 2 l + 1}) _{j - l} }}\\
\times
{} _{4} F _{3} \!\left [ \matrix { {{-j}\over 2} + {l\over 2}, {1\over 2}
       - {j\over 2} + {l\over 2}, -i - 2 j - x, 1 + 2 j + 2 l + 2 x}\\ { 1
       - j + l, {1\over 2} - {i\over 2} + l + {x\over 2}, 1 - {i\over 2} + l +
       {x\over 2}}\endmatrix ; {\displaystyle 1}\right ] \\
- 
  {{       ({ \textstyle x - i + 2 j + 1}) _{n+i} \,
      ({ \textstyle x + i + 2 l + 1}) _{-i + j - l + n} }\over 
    {({ \textstyle 2 x + 2 j + 2 l + 1}) _{j - l} }}\hskip6cm\\
\times
{} _{4} F _{3} \!\left [ \matrix { {{-j}\over 2} + {l\over 2}, {1\over 2}
       - {j\over 2} + {l\over 2}, i - 2 j - x, 1 + 2 j + 2 l + 2 x}\\ { 1
       - j + l, {1\over 2} + {i\over 2} + l + {x\over 2}, 1 + {i\over 2} + l +
       {x\over 2}}\endmatrix ; {\displaystyle 1}\right ]\\
 +   ({ \textstyle x - i + 2 j + 1}) _{n+i} \,
   ({ \textstyle x + i + j + l + 1}) _{n-i}\hskip4cm  \\
- ({ \textstyle x + i + 2 j + 1}) _{n-i} \,
     ({ \textstyle x - i + j + l + 1}) _{n+i}   =0.
\endmultline\tag\ZF$$
Both $_4F_3$-series can be summed by means of a $_4F_3$-summation
which appears in a paper by Andrews and Burge \cite{\AnBuAA,
Lemma~1} (see \cite{\KratBD, Lemma~A3} for a simpler proof),
$${}_4F_3\!\[\matrix -\frac {N} {2},\frac {1} {2}-\frac {N} {2},-a,a+b\\
1-N,\frac {b} {2},\frac {1} {2}+\frac {b} {2}\endmatrix; 1\]=\frac {(a+b)_N} {(b)_N}+
\frac {(-a)_N} {(b)_N},$$
where $N$ is a positive integer. We have to apply the case where
$N=j-l$. This is indeed a positive integer because of our
assumption $j>l$. Some simplification then leads to (\ZF). 

This shows that
$\prod _{j=1} ^{n}(x-y+2j+1)_j$ divides $D_n(x,y)$.

\smallskip
{\it Step 2. $\prod _{j=1} ^{n}(x+y+j+1)_j$ is a factor of
$D_n(x,y)$}. Let us concentrate on a typical factor $(x+y+j+l)$,
$1\le j\le n$, $1\le l\le j$. We claim that for each
such factor there is a linear combination of the columns that vanishes
if the factor vanishes. More precisely, we claim that for any $j,l$
with $1\le j\le n$, $1\le l\le j$ there holds
$$\multline 
\sum_{s = 1 + j - l}^{j}
      \(-\frac {1}{4}\)^{j-s}{{l-1}\choose {s+l - j-1}}\\
\cdot{{        
       ({ \textstyle x + n + 2 s + 1}) _{2 j - 2 s} \,
       ({ \textstyle 2 x + 3 j + l + s + 1}) _{j - s} }\over 
     {({ \textstyle x + j + s + {1\over 2}}) _{j - s} \,
       ({ \textstyle x+j + l + s }) _{j - s} \,
       ({ \textstyle x+j + l - n + s }) _{j - s} }}\\
\cdot(\text {column
$s$ of $D_n(x,-x-j-l)$}) =0.
\endmultline\tag\ZG$$
In order to verify (\ZG), we have to check
$$\multline 
\sum_{s = 1 + j - l}^{j}
      \(-\frac {1}{4}\)^{j-s}{{l-1}\choose {s+ l-j-1}}\\
\cdot{{       
       ({ \textstyle x + n + 2 s + 1}) _{2 j - 2 s} \,
       ({ \textstyle 2 x + 3 j + l + s + 1}) _{j - s} }\over 
     {({ \textstyle x+ j + s + {1\over 2}}) _{j - s} \,
       ({ \textstyle x+j + l + s }) _{j - s} \,
       ({ \textstyle x+j + l - n + s }) _{j - s} }}\\
\cdot    \big( ({ \textstyle -x + i - j - l - s +1}) _{n-i} \,
          ({ \textstyle x - i + 2 s + 1}) _{n+i} \\
 -      ({ \textstyle -x - i - j - l - s +1}) _{n+i} \,
          ({ \textstyle x + i + 2 s + 1}) _{n-i}  \big)=0,
\endmultline$$
which is (\ZG) restricted to the $i$-th row.
Equivalently, using hypergeometric notation,
this means to check
$$\multline 
 {{\left( -1 \right) }^l} 
  {{        ({ \textstyle -i - 2 j - x}) _{n+i} \,
      ({ \textstyle 3 + i + 2 j - 2 l + x}) _{-2 - i + 2 l + n} \,
      ({ \textstyle 2 + 4 j + 2 x}) _{l-1} }\over 
    {{4^{l-1}} ({ \textstyle 1 + 2 j + x}) _{l-1} \,
      ({ \textstyle {3\over 2} + 2 j - l + x}) _{l-1} \,
      ({ \textstyle 1 + 2 j - n + x}) _{l-1} }}\hskip2cm\\
\times
{} _{4} F _{3} \!\left [ \matrix { 1 - l, {3\over 2} + 2 j - l + x, 1 +
       2 j + x, 1 + i + 2 j + x}\\ { 2 + 4 j + 2 x, 2 + {i\over 2} + j - l
       + {x\over 2}, {3\over 2} + {i\over 2} + j - l + {x\over 2}}\endmatrix ;
       {\displaystyle 1}\right ]\\
- {{\left( -1 \right) }^l} 
{{       ({ \textstyle i - 2 j - x}) _{n-i} \,
      ({ \textstyle 3 - i + 2 j - 2 l + x}) _{-2 + i + 2 l + n} \,
      ({ \textstyle 2 + 4 j + 2 x}) _{l-1} }\over 
    {{4^{l-1}} ({ \textstyle 1 + 2 j + x}) _{l-1} \,
      ({ \textstyle {3\over 2} + 2 j - l + x}) _{l-1} \,
      ({ \textstyle 1 + 2 j - n + x}) _{l-1} }}\\
\times
{} _{4} F _{3} \!\left [ \matrix { {3\over 2} + 2 j - l + x, 1 + 2 j +
       x, 1 - i + 2 j + x, 1 - l}\\ { {3\over 2} - {i\over 2} + j - l +
       {x\over 2}, 2 - {i\over 2} + j - l + {x\over 2}, 2 + 4 j +
       2 x}\endmatrix ; {\displaystyle 1}\right ]=0.
\endmultline\tag\ZH$$
In order to establish (\ZH) we apply Bailey's transformation for 
balanced $_4F_3$-series
(see \cite{\SlatAC, (4.3.5.1)}),
$$\multline
{} _{4} F _{3} \!\left [ \matrix { a, b, c, -N}\\ { e, f, 1 + a + b + c - e -
   f - N}\endmatrix ; {\displaystyle 1}\right ]  \\=
{{( e-a)_N\, (f-a) _{N}}\over  {(e)_N\,( f) _{N}}}
  {} _{4} F _{3} \!\left [ \matrix { -N, a, 1 + a + c - e - f - N, 1 + a + b -
    e - f - N}\\ { 1 + a + b + c - e - f - N, 1 + a - e - N, 1 + a - f -
    N}\endmatrix ; {\displaystyle 1}\right ]   ,
\endmultline$$
where $N$ is a nonnegative integer, to the second $_4F_3$-series in
(\ZH). Thus it is converted into the first $_4F_3$-series, and it is
routine to check that also the remaining terms that go with the
$_4F_3$-series agree. So, the two terms on the left hand side of (\ZH)
cancel each other, as desired. 

This establishes that $\prod _{j=1} ^{n}(x+y+j+1)_j$ divides
$D_n(x,y)$.

\smallskip
{\it Step 3. $\prod _{i=1} ^{n}(x+2y+3i+1)_{n-i}$ is a factor of
$D_n(x,y)$}. This is the most difficult part of the proof of (\ZC).
Trials of finding linear combinations of columns that vanish resulted in
extremely messy expressions. So, we decided to work with linear
combinations of rows this time. Still, the coefficients are not as
``nice" as in Steps~1 and 2.

Let us concentrate on a typical factor $(x+2y+3i+l)$,
$1\le i\le n$, $1\le l\le n-i$. We claim that for each
such factor there is a linear combination of the rows that vanishes
if the factor vanishes. More precisely, we claim that for any $i,l$
with $1\le i\le n$, $1\le l\le n-i$ there holds
$$\sum _{k=1} ^{i+l}\frac {(k+i+l+1)_{i+l-k}} {(i+l-k)!}\,P_l(2i,i+l-k)
\cdot(\text {row $k$ of $D_n(-2y-3i-l,y)$}) =0,
\tag\ZI$$
where $P_l(e,f)$ is the polynomial
$$P_l(e,f)=\sum _{r=0} ^{2l+1}a_r \,(e)_r \,(-f)_{2l+1-r},\tag\ZJ$$
with the expansion coefficients $a_r$ given by
$$a_r=\langle x^r\rangle \,\big((x^2+x+1)^{l-1}(2x+1)(x+2)(x-1)\big).\tag\ZK$$
Here, $\langle x^r\rangle g(x)$ denotes the coefficient of $x^r$ in
$g(x)$.

By specializing (\ZI) to the $j$-th column, splitting the resulting
sum into two parts in the obvious way, and then moving one sum to
the right hand side, we see that in order to verify (\ZI), we have to check
$$\multline 
\sum _{k=1} ^{i+l}\frac {(k+i+l+1)_{i+l-k}} {(i+l-k)!}\,P_l(2i,i+l-k)\,
(-2y-3i-l+2j-k+1)_{n+k}\,(y+k-j+1)_{n-k}\\
=\sum _{k=1} ^{i+l}\frac {(k+i+l+1)_{i+l-k}}
{(i+l-k)!}\,P_l(2i,i+l-k)\,
(-2y-3i-l+2j+k+1)_{n-k}\,(y-k-j+1)_{n+k},
\endmultline$$
or, after adding one more term as first 
summand on both sides, equivalently,
$$\multline 
\sum _{k=0} ^{i+l}\frac {(k+i+l+1)_{i+l-k}} {(i+l-k)!}\,P_l(2i,i+l-k)\,
(-2y-3i-l+2j-k+1)_{n+k}\,(y+k-j+1)_{n-k}\\
=\sum _{k=0} ^{i+l}\frac {(k+i+l+1)_{i+l-k}}
{(i+l-k)!}\,P_l(2i,i+l-k)\,
(-2y-3i-l+2j+k+1)_{n-k}\,(y-k-j+1)_{n+k}.
\endmultline\tag\ZL$$
Empirically, we discovered that apparently both sums in
(\ZL) are indefinitely summable (``Gosper-summable"; see
\cite{\GrKPAA, \S5.7}, \cite{\PeWZAA, \S II.5}). It is exactly this fact
which makes (\ZL) tractable. 

In the following we will show that the sums in (\ZL) are equal,
however, without exhibiting an explicit expression for the sums.
Instead, what we will do is to read through Gosper's algorithm
\cite{\GospAB} (see also \cite{\GrKPAA, \S5.7}, \cite{\PeWZAA, \S II.5}), which is
an algorithm that solves the problem of indefinite summation for
hypergeometric sums. (For any fixed $l$, our sums in (\ZL) belong to the
category of hypergeometric sums.) In the course of reading through
Gosper's algorithm it will emerge that the sums on both sides of
(\ZL) must be equal.

Let us recall what Gosper's algorithm does and how it works. Let $t(k)$
be a ``hypergeometric term", i.e., be a term such that the ratio
$t(k+1)/t(k)$ is a rational function in $k$. Then the Gosper
algorithm will find a hypergeometric term $T(k)$ (if it exists)
satisfying
$$t(k)=T(k+1)-T(k).\tag\ZM$$
The upshot of this is that then the indefinite summation of the term
$t(k)$ can be easily carried out,
$$\sum _{k=A} ^{B}t(k)=T(B+1)-T(A).\tag\ZN$$
The term $T(k)$ is found in the following way. First, one finds
polynomials $p(k)$, $q(k)$, and $r(k)$ such that
$$\frac {t(k+1)} {t(k)}=\frac {p(k+1)} {p(k)}\frac {q(k)}
{r(k+1)},\tag\ZO$$
where $q(k)$ and $r(k)$ have the property that whenever $(k+\al)\mid
q(k)$ and $(k+\be)\mid r(k)$ then the difference $\al-\be$ must not
be a positive integer. Next, one finds a polynomial $s(k)$
satisfying the recurrence relation
$$p(k)=q(k)s(k+1)-r(k)s(k)\tag\ZP$$
for all $k$. The term $T(k)$ is then given by
$$T(k)=\frac {r(k)\,s(k)} {p(k)}t(k).\tag\ZR$$

Now let us carry out this program with the summands in (\ZL). First,
let $t(k)=t_1(k)$, where $t_1(k)$ is the summand of the sum on the
left hand side of (\ZL),
$$t_1(k)=\frac {(k+i+l+1)_{i+l-k}} {(i+l-k)!}\,P_l(2i,i+l-k)\,
(-2y-3i-l+2j-k+1)_{n+k}\,(y+k-j+1)_{n-k}.$$
Then (\ZO) holds with $p(k)=p_1(k)$, $q(k)=q_1(k)$, $r(k)=r_1(k)$,
where $p_1(k)=P_l(2i,i+l-k)$,
$q_1(k)=(i+l-k)(-2y-3i-l+2j-k)$, and $r_1(k)=(i+l+k)(y-j+k)$. So, next we
have to find a polynomial $s_1(k)$ satisfying the recurrence
$$P_l(2i,i+l-k)=(i+l-k)(-2y-3i-l+2j-k)s_1(k+1)-(i+l+k)(y-j+k)s_1(k).\tag\ZS$$
For each specific instance of $i$ and $l$ this is just routine.
However, we were not able to find an explicit formula for $s_1(k)$ in
general. Fortunately, we do not need such an
explicit expression. Assuming that we have found a
polynomial $s_1(k)$ satisfying (\ZS), by (\ZN) and (\ZR) we have
$$\align \sum _{k=0} ^{i+l}\frac {(k+i+l+1)_{i+l-k}} {(i+l-k)!}&\,P_l(2i,i+l-k)\,
(-2y-3i-l+2j-k+1)_{n+k}\,(y+k-j+1)_{n-k}\\
&=
\frac {r_1(i+l+1)\,s_1(i+l+1)} {p_1(i+l+1)}t_1(i+l+1)-\frac {r_1(0)\,s_1(0)}
{p_1(0)}t_1(0)\\
&=-\frac {(i+l)_{i+l+1}}
{(i+l)!}(-2y-3i-l+2j+1)_n\,(y-j)_{n+1}\,s_1(0),\tag\ZT
\endalign$$
the last line being due to the fact that $t_1(i+l+1)=0$.

On the other hand, for $t(k)=t_2(k)$, where $t_2(k)$ is the summand of the sum on the
right hand side of (\ZL),
$$t_2(k)=\frac {(k+i+l+1)_{i+l-k}}
{(i+l-k)!}\,P_l(2i,i+l-k)\,
(-2y-3i-l+2j+k+1)_{n-k}\,(y-k-j+1)_{n+k}$$
we may choose $p(k)=p_2(k)$, $q(k)=q_2(k)$, $r(k)=r_2(k)$,
where $p_2(k)=P_l(2i,i+l-k)$,
$q_2(k)=(i+l-k)(y-j-k)$, and $r_2(k)=(i+l+k)(-2y-3i-l+2j+k)$.
So, here we have to find a polynomial $s_2(k)$ satisfying the recurrence
$$P_l(2i,i+l-k)=(i+l-k)(y-j-k)s_2(k+1)-(i+l+k)(-2y-3i-l+2j+k)s_2(k).\tag\ZU$$
Again, this is just routine for each specific instance of $i$ and $l$,
but we do not know an explicit formula for $s_2(k)$ in general. 
Assuming that we have found a
polynomial $s_2(k)$ satisfying (\ZU), by (\ZN) and (\ZR) we have
$$\align \sum _{k=0} ^{i+l}\frac {(k+i+l+1)_{i+l-k}} {(i+l-k)!}&\,P_l(2i,i+l-k)\,
(-2y-3i-l+2j+k+1)_{n-k}\,(y-k-j+1)_{n+k}\\
&=
\frac {r_2(i+l+1)\,s_2(i+l+1)} {p_2(i+l+1)}t_2(i+l+1)-\frac {r_2(0)\,s_2(0)}
{p_2(0)}t_2(0)\\
&=-\frac {(i+l)_{i+l+1}}
{(i+l)!}(-2y-3i-l+2j)_{n+1}\,(y-j+1)_{n}\,s_2(0),\tag\ZV
\endalign$$
the last line being due to the fact that also $t_2(i+l+1)=0$.

In order to relate $s_2(k)$ to $s_1(k)$, we make the following
observation: We set $s_2(k)=\tilde s_2(-k+1)$, substitute this in the
recurrence (\ZU), then replace $k$ by $-k$ and change the sign on
both sides of (\ZU). Thus we obtain for
$\tilde s_2(k)$ the recurrence
$$-P_l(2i,i+l+k)=(i+l-k)(-2y-3i-l+2j-k)\tilde
s_2(k+1)-(i+l+k)(y-j+k)\tilde s_2(k).\tag\ZW$$
This is almost the same recurrence as the recurrence (\ZS) for
$s_1(k)$! It is only the term on the left hand side which is
different! But, in fact, there is no difference: We claim that:
\roster
\item"\hskip1cm Claim 1:\hss\hss" \kern1.3cm We have $P_l(e,e+2l-f)=-P_l(e,f)$.
\item"\hskip1cm Claim 2:\hss\hss" \kern1.3cm There exists a unique solution for the recurrence (\ZS).
\endroster
Let us for the moment assume that these claims have been already
established. Then, because of Claim~1, the recurrences (\ZS) and (\ZW)
are indeed the same. Furthermore, thanks to Claim~2, there does exist
a unique solution for the recurrence (\ZS), and so also for (\ZW). Hence,
the solutions must be the same, i.e., $s_1(k)=\tilde s_2(k)$, which
means $s_1(k)=s_2(1-k)$. In particular, we have $s_1(1)=s_2(0)$. 
A further fact, which follows immediately from Claim~1 on replacing
$e$ by $2e$ and setting
$f=e+l$, is that $P_l(2e,e+l)=0$. Therefore,
by setting $k=0$ in (\ZS), we obtain
$$0=(i+l)(-2y-3i-l+2j)s_1(1)-(i+l)(y-j)s_1(0).$$
{}From this equation, and the previous observation that $s_1(1)=s_2(0)$, we infer 
$$s_1(0)=\frac {(-2y-3i-l+2j)} {(y-j)}s_1(1)=\frac {(-2y-3i-l+2j)}
{(y-j)}s_2(0).$$
Substitution of this relation in (\ZT) gives
$$\align \sum _{k=0} ^{i+l}\frac {(k+i+l+1)_{i+l-k}} {(i+l-k)!}&\,P_l(2i,i+l-k)\,
(-2y-3i-l+2j-k+1)_{n+k}\,(y+k-j+1)_{n-k}\\
&=-\frac {(i+l)_{i+l+1}}
{(i+l)!}(-2y-3i-l+2j)_{n+1}\,(y-j+1)_{n}\,s_2(0).
\endalign$$
Comparison of this identity with (\ZV) shows that indeed the sums on
both sides of (\ZL) are identical. This would prove (\ZL).

So it remains to settle Claims~1 and 2. 

We begin with Claim~1. By the definition (\ZJ) of $P_l(e,f)$, we have
$$\align P_l(e,e+2l-f)
&=\sum _{r=0} ^{2l+1}a_r \,(e)_r \,(-e-2l+f)_{2l+1-r}\\
&=\sum _{r=0} ^{2l+1}a_r \,(e)_r \,(-1)^{r+1}\,(e+r-f)_{2l+1-r},
\endalign$$
where the coefficients $a_r$ are given by (\ZK). Next we use the
Chu--Vandermonde summation (see e.g\. \cite{\GrKPAA, 
\S5.1, (5.27)}) in the form
$$\sum _{s=0} ^{N}\binom {N}s(x)_s\,(y)_{N-s}=(x+y)_N,$$
with $N=2l+1-r$, $x=e+r$, and $y=-f$.
Thus,
$$\align P_l(e,e+2l-f)
&=\sum _{r=0} ^{2l+1}a_r \,(e)_r \,(-1)^{r+1}\,\sum _{s=0}
^{2l+1-r}\binom {2l+1-r}s(e+r)_s\,(-f)_{2l+1-r-s}\\
&=-\sum _{m=0} ^{2l+1}(e)_m\,(-f)_{2l+1-m}\sum _{r=0} ^{m}\binom
{2l+1-r}{m-r}(-1)^ra_r.
\endalign$$
Therefore, Claim~1 will follow immediately, if we are able to show that
$$\sum _{r=0} ^{m}\binom
{2l+1-r}{m-r}(-1)^ra_r=a_m.\tag\ZX$$
This can be readily done by using generating functions. The
definition (\ZK) of the coefficients $a_r$ is equivalent to
$$\sum _{r=0} ^{\infty}a_rx^r=(x^2+x+1)^{l-1}(2x+1)(x+2)(x-1).\tag\ZXa$$
Let us denote the right hand side of this equation by $A(x)$. 
Now we multiply both
sides of (\ZX) by $x^m$, and we sum over all $m=0,1,\dots$ We obtain
$$\sum _{m=0} ^{\infty}\sum _{r=0} ^{m}\binom
{2l+1-r}{m-r}(-1)^ra_r=A(x),$$
and after interchanging summations on the left hand side and summing
the (now) inner sum over $m$ by means of the binomial theorem,
$$(1+x)^{2l+1}A\left(-\frac {x} {1+x}\right)=A(x).$$
It is trivial to verify this equation. Thus also the equivalent
equation (\ZX) must be true.
Due to the preceding considerations, 
this completes the proof of Claim~1.

Next we turn to Claim~2. We show that there is a unique polynomial $s_1(k)$
of degree $2l$ that satisfies the recurrence (\ZS). (We leave it as
an exercise that the ``degree calculus" of the Gosper algorithm shows
that if there is a solution to the recurrence (\ZS) then it has to be
a polynomial of degree at most $2l$.) So, let
$s_1(k)=\sum _{m=0} ^{2l}c(m) (k-i-l)_m$. We substitute this into
(\ZS), then expand everything with respect to the basis $(k-i-l)_m$,
$m=0,1,\dots$ (for the space of polynomials in $k$), and finally
compare coefficients of $(k-i-l)_m$ on both sides of (\ZS). This
leads to the following system of equations for the coefficients $c(m)$:
$$\multline
a_{2l+1-m}(2i)_{2l+1-m}=(y+i-l-j+m-1)\,c(m-1)-(2i+2l-m)(y+i+l-j-m)\,c(m),\\
m=0,1,\dots,2l+1,
\endmultline\tag\ZY$$
where, by convention, we put $c(-1)=c(2l+1)=0$. For convenience, we
set 
$$c(m)=(2i)_{2l-m}\,(y+i-l-j)_{2l-m}\,(y+i-l-j)_m \,\tilde c(m).$$
By substituting this in (\ZY), we obtain the simpler system of
equations
$$
\frac {a_m} {(y+i-l-j)_{2l-m+1}\,(y+i-l-j)_m}=\tilde c(m-1)-\tilde 
c(m),\quad \quad  m=0,1,\dots,2l+1.
\tag\ZZ$$
This is a system of $2l+2$ equations for $2l+1$ variables. 
(Recall the convention $c(-1)=c(2l+1)=0$, which of course implies
$\tilde c(-1)=\tilde c(2l+1)=0$.) So, it is
overdetermined. It is easy to see that this inhomogeneous system of linear
equations has a (unique) solution if and only if the sum of the left hand
sides of (\ZZ) over all $m$ equals 0, i.e., if and only if
$$\sum _{m=0} ^{2l+1}\frac {a_m}
{{(y+i-l-j)_{2l-m+1}\,(y+i-l-j)_m}}=0.\tag\Za$$
This would follow immediately from the antisymmetry property
$a_m=-a_{2l+1-m}$, because then the $m$-th and $(2l+1-m)$-th summand
in the sum in (\Za) would cancel each other. Indeed, the substitution
$x\to 1/x$ in (\ZXa) yields $a_m=-a_{2l+1-m}$.
Therefore, the system of equations (\ZY) has indeed a unique
solution, which implies that there is a unique polynomial $s_1(k)$
satisfying the recurrence (\ZS), which is exactly the assertion of
Claim~2.

\smallskip
The proof that $\prod _{i=1} ^{n}(x+2y+3i+1)_{n-i}$ divides
$D_n(x,y)$ is now complete.

\smallskip
{\it Step 4. $D_n(x,y)$ is a polynomial in $x$ of maximal degree
$n(3n+1)/2$, and the same is true for the maximal degree of $D_n(x,y)$
in $y$}. This is because each term in the defining expansion
of the determinant $D_n(x,y)$ has degree $n(3n+1)/2$ in $x$, and the
same in $y$. Since the right hand side of (\ZC), which
by Steps~1--3 divides $D_n(x,y)$ as a polynomial in $x$ and $y$,
also has degree $n(3n+1)/2$ in $x$, respectively $y$, $D_n(x,y)$ and the right hand side of
(\ZC) differ only by a multiplicative constant.

\smallskip
{\it Step 5. 
The evaluation of the multiplicative constant.} By the preceding
steps we know that (\ZD) holds. In particular, if we set $y=0$, we
have
$$\det_{1\le i,j\le n}\((x+2j-i+1)_{n+i}\,(i-j+1)_{n-i}\)
=C(n)\prod _{j=1} ^{n}{(x+j+1)_{n+j}} .\tag\Zb$$
(The reader should be aware that the second term in the determinant
$D_n(x,y)$, as given by (\ZC), vanishes for $y=0$ because of the
presence of the factor $(y-i-j+1)_{n+i}$.) 
The determinant on the left hand side of (\Zb) is a lower triangular
matrix, hence it equals the product of its diagonal entries, which is
$\prod _{j=1} ^{n}(x+j+1)_{n+j}\,(n-j)!$. Therefore
$C(n)$ is equal to $\prod _{j=1} ^{n}(n-j)!=\prod _{j=1} ^{n}(j-1)!$.

\medskip
This finishes the proof of (\ZC) and thus of the Theorem.\quad \quad \qed
\enddemo

\remark{Acknowledgement} We are grateful to the referee for a
simplification of our original proof of (\ZX).
\endremark

\enddocument